\documentclass[11pt]{article}
\usepackage{theorem,amsmath,amssymb,xypic}
\xyoption{curve}
\theoremstyle{change}
\theorembodyfont{\itshape}
\newtheorem{thm}[subsection]{Theorem}
\newtheorem{prop}[subsection]{Proposition}
\newtheorem{lemma}[subsection]{Lemma}

\newtheorem{defn}[subsection]{Definition}
{\theorembodyfont{\rmfamily}

}

\makeatletter
\renewcommand{\subsection}{\@startsection%
{subsection}{2}{0mm}{\baselineskip}{-1em}%
{\normalfont\normalsize\bfseries}}
\renewcommand{\subsubsection}{\@startsection%
{subsubsection}{3}{0mm}{\baselineskip}{-1em}%
{\normalfont\normalsize\textit}}
\makeatother

\numberwithin{equation}{subsection}

\renewcommand{\labelenumi}%
  {\stepcounter{equation}\rm(\theequation)}

\newcommand{\demobox}{\vrule height6pt width6pt depth0pt}

\newenvironment{demo}{\noindent{\it Proof.}}
{{\unskip\nobreak\hfil\qquad
\demobox\parfillskip=0pt\par}
\medskip}

\newcommand{\Mod}{{\operatorname{Mod}}}

\newcommand{\tens}{\otimes}

\newcommand{\bF}{{\mathbb{F}}}

\newcommand{\bC}{{\mathbb{C}}}

\newcommand{\bQ}{{\mathbb Q}}

\newcommand{\bZ}{{\mathbb Z}}

\newcommand{\Xto}{\xrightarrow}

\newcommand{\isom}{\Xto{\sim}}
\newcommand{\sset}{\subseteq}

\newcommand{\perf}{\mathit{perf}}

\newcommand{\car}{\mathrm{char}}
\newcommand{\rig}{\mathit{rig}}
\newcommand{\Fet}{\mathrm{Fet}}
\newcommand{\fet}{\mathrm{fet}}

\newcommand{\tor}{\mathit{tor}}
\newcommand{\fTd}{\mathit{ftd}}
\newcommand{\et}{\mathrm{et}}
\newcommand{\cohdim}{\mathrm{coh.dim}}
\newcommand{\crys}{\mathit{crys}}
\newcommand{\cryset}{\mathit{crys-et}}
\newcommand{\cG}{\mathcal{G}}

\begin{document}

\title{\bf Integral $p$-adic cohomology theories}
\author{Tomoyuki Abe and Richard Crew}
\date{}

\maketitle

\section{Introduction}
\label{sec:intro}

Suppose $\Lambda$ is a complete local ring with residue field $k$. Can
we expect that there is a ``reasonable'' cohomology theory
$H^\cdot(X,\Lambda)$ with its values in $\Lambda$-modules for
separated schemes of finite type over a field $k$? Here ``reasonable''
means, in a first approximation that it has all the usual properties
of integral Betti cohomology when $k=\bC$. For example if
$\Lambda=\bZ_\ell$ with $\ell$ different from $p=\car(k)$, $\ell$-adic
cohomology $H^\cdot_\et(X,\bZ_\ell)$ is such a theory. If $\Lambda=W$
is a Cohen ring of $k$, crystalline cohomology gives such a theory
when $X/k$ is proper and smooth, but not in more general situations:
Berthelot found that the torsion of $H^1_\crys(X,W)$ could be infinite
if $X$ was even mildly singular. Of course if we replace $W$ by its
fraction field $K$, rigid cohomology has the desired properties.

In the article in which he first discussed such theories Grothendieck
emphasized the importance of $p$-torsion phenomena, which is only
visible in an ``integral'' $p$-adic theory such as crystalline
cohomology:
\begin{quotation}
  Such a theory should associate to each scheme $X$ of finite type
  over a perfect field $k$ of characteristic $p>0$, cohomology groups
  which are modules over an integral domain, whose quotient field is
  of \textit{characteristic $0$}, and which satisfy all the desirable
  formal properties (functoriality, finite-dimensionality...). This
  cohomology should also, most importantly, explain torsion phenomena,
  and in particular \textit{$p$-torsion.} (\cite[\S1.7]{Grothendieck},
  emphasis in the original).
\end{quotation}
The work of Illusie and others on the de Rham-Witt complex shows the
immense richness of $p$-torsion phenonmena, but of course only in the
proper smooth case.

When $k$ is perfect and $X/k$ is smooth Davis, Langer and Zink
\cite{DLZ} have defined an overconvergent version of the de Rham-Witt
complex $W^\dagger\Omega^\cdot_{X/W(k)}$ and constructed an isomorphism
\begin{displaymath}
  H^\cdot(X,W^\dagger\Omega^\cdot_{X/W(k)})\tens\bQ\simeq H^\cdot_\rig(X)
\end{displaymath}
when $X$ is quasiprojective. It's evidently too much to hope that the
$W(k)$-modules $H^\cdot(X,W^\dagger\Omega^\cdot_{X/W(k)})$ are
finitely generated, and we will see that this is \textit{never} true
if $X$ is a smooth affine curve. Davis, Langer and Zink made the more
reasonable conjecture that the image of
$H^\cdot(X,W^\dagger\Omega^\cdot_{X/W(k)})$ in $H^\cdot_\rig(X)$ is
finitely generated as a $W(k)$-module, but in a recent preprint Ertl
and Shiho \cite{Ertl-Shiho} have produced counterexamples to this
assertion as well. In any case one loses all torsion information by
replacing $H^\cdot(X,W^\dagger\Omega^\cdot_{X/W(k)})$ by its image in
$H^\cdot_\rig(X)$.

The main result of this paper is a nonexistence theorem for certain
theories of this sort. We assume that a theory has a comparison
theorem with rigid cohomology (or rigid cohomology with compact
support), that the cohomology groups are finitely generated
$\Lambda$-modules, and -- a natural but important additional
restriction -- the theory is compatible with finite \'etale descent in
a manner to be explained later. The result is that there is no such
theory even for affine curves. We do not rule out the possibility of
theories satisfying weaker descent conditions, such as cdh-descent.
In fact in another recent preprint \cite{Ertl-Shiho-Sprang} Ertl,
Shiho and Sprang construct a ``good'' integral $p$-adic theory under
certain assumptions concerning resolution of singularities in positive
characteristic. This construction uses cdh-descent and so does not
contradict our result if their hypotheses on resolution hold. They
also consider a theory based on simplicial generically \'etale
hypercovers, and show that it is \textit{not} independent of the
choice of hypercover, and thus does not provide a theory compatible
with finite \'etale descent.  We should remark finally that Bhargav
Bhatt has also given a proof of the nonexistence of this sort of
theory (private communication) by considering Artin-Schreier covers of
the the affine line by itself; he shows that the existence of such a
theory would imply that the affine line has Euler characteristic
$0$. The method of the present paper on the other hand shows that the
misbehavior of $p$-adic cohomology theories is ubiquitous, at least as
far as curves are concerned; see the example at the end of section
\ref{sec:examples}.

The original version of this result was inspired by discussions at the
conference ``$p$-adic cohomology and arithmetic applications'' at BIRS
in October 2017. Later versions were produced during visits of the
second author to the Kavli IPMU in November 2019 and the University of
Rennes I in May 2022, and he would like to thank both instutiions for
their hospitality. The first author is supported by JSPS KAKENHI Grant
Numbers 16H05993, 18H03667, 20H01790.

\section{Globally Perfect Models}
\label{sec:models}

In what follows $\Lambda$ is a complete noetherian local ring with
residue field $k$ and fraction field $K$.\ As usual $D_\perf(\Lambda)$
is the triangulated category of perfect complexes of
$\Lambda$-modules, which since $\Lambda$ is local means that an object
of $D_\perf(\Lambda)$ is quasi-isomorphic to a bounded complex of
\textit{free} $\Lambda$-modules.

We have in mind the following requirements on a cohomology theory
$H^\cdot(X)$ on a subcategory of the category of $k$-schemes of finite
type:
\begin{itemize}
\item $H^\cdot(X)$ can be used to compute rigid cohomology (or rigid
  cohomology with compact supports) via a suitable comparison
  theorem.
\item $H^\cdot(X)$ may be computed as the cohomology of an object of
  $D_\perf(\Lambda)$; 
\item The $H^\cdot(X)$ are compatible with finite \'etale descent, in a
  sense to be explained presently.
\end{itemize}

When $\Lambda$ is regular the second requirement reduces to the
condition that the $H^n(X)$ be finitely generated $\Lambda$-modules,
and vanish for $|n|\gg0$.

To explain the last condition we use the \textit{finite \'etale site},
referring to the book of Abbes, Gros and Tsuji \cite[Ch. VI \S9]{AGT}
for proofs of the assertions used below. If $X$ is a scheme the site
$\Fet(X)$ is the category of finite \'etale morphisms $Y\to X$, and
the coverings are surjective morphisms. The associated topos will be
written $X_\fet$. A morphism $\pi:Y\to X$ induces a morphism
$\pi_\fet:Y_\fet\to X_\fet$ of topoi; in what follows we will
abbreviate the associated functors $\pi_\fet^*$ and $\pi_{\fet *}$ by
$\pi^*$ and $\pi_*$ (in other words the latter refer to the morphism
of finite \'etale topoi, not \'etale topoi).  There is also a
projection $\rho_X:X_\et\to X_\fet$ such that
$\pi_\fet\rho_Y\simeq\rho_X\pi_\et$ in $\Fet(X)$ for any $\pi:Y\to
X$. If $X$ is a coherent scheme with finitely many components (e.g. if
$X$ is of finite type over a field) the inverse image
$\rho_X^*:X_\fet\to X_\et$ is fully faithful
(\cite[Prop. VI.9.18]{AGT}).

\begin{defn}\label{defn:globally-perfect}
  An object $M$ of $D^+(X_\fet,\Lambda)$ is \emph{globally perfect} if
  there are integers $a\le b$ such that for all $\pi:Y\to X$ in
  $\Fet(X)$, $R\Gamma(Y,\pi^*M)$ is an object of
  $D^{[a,b]}_\perf(\Lambda)$.
\end{defn}

We denote by $\Mod^\cdot_K$ the category of $\bZ$-graded $K$-vector
spaces. The next definition formulates our notion of what it means for
a cohomology theory with values in $K$-vector spaces to have a
``good'' integral model compatible with finite \'etale descent.

\begin{defn}
  Suppose $X$ is a $k$-scheme of finite type and
  $H^\cdot:\Fet(X)\to\Mod^._K$ is a functor. A \emph{globally perfect
    model of $H^\cdot$} is a globally perfect $M$ in
  $D^+(X_\fet,\Lambda)$ and a functorial isomorphism
  \begin{displaymath}
    H^\cdot(Y)\isom H^\cdot(R\Gamma(Y,\pi^*M))\tens_\Lambda K
  \end{displaymath}
  for all $Y\to X$ in $\Fet(X)$.
\end{defn}

Note that if $\Lambda$ is regular the condition on
$R\Gamma(Y,\pi^*M)$ is equivalent to saying that the $H^p(Y,\pi^*M)$
are finitely generated and vanish outside of a finite range depending
only on $X$.

\begin{thm}\label{thm:perfection}
  Suppose $\Lambda$ is noetherian, $X$ is \emph{affine}, $k_\et$ has
  finite cohomological dimension, $M$ is a globally perfect object of
  $D(X_\fet,\Lambda)$ and $\pi:Y\to X$ is finite \'etale Galois with
  group $G$. There is an object $M_Y$ of $D_\perf(\Lambda[G])$ whose
  image under the forgetful functor
  $D_\perf(\Lambda[G])\to D_\perf(\Lambda)$ is isomorphic to
  $R\Gamma(Y,\pi^*M)$, and $H^\cdot(M_Y)\simeq H^\cdot(Y,\pi^*M)$ as
  $\Lambda[G]$-modules.
\end{thm}
We could express the conclusion of the theorem by saying that
``$R\Gamma(Y,\pi^*M)$ with its $G$-action is a perfect complex of
$\Lambda[G]$-modules.''

For the proof of theorem \ref{thm:perfection} we will need some facts
about cohomological dimension. Recall that if $R$ is a ring and $X$ is
a topos, the $R$-cohomological dimension of $X$ is the smallest
integer $d$ such that $H^n(X,F)=0$ for all $n>d$ and every $R$-module
$F$ in $X$. The $\bZ$-cohomological dimension of $X$ will also be
called simply the cohomological dimension and will be denoted by
$\cohdim(X)$.

\begin{lemma}\label{lemma:coh-dim}
  Suppose $k$ is a field of characteristic $p>0$ such that $k_\et$ has
  finite cohomological dimension. Let $X$ be a $k$-scheme of finite
  type.
  \begin{itemize}
  \item The cohomological dimension of $X_\et$ is
    finite.
  \item If $X$ is affine the cohomological dimension of $X_\fet$
    is finite.
  \end{itemize}
\end{lemma}
\begin{demo}
  If $k$ is separably closed, the first assertion is a theorem of
  Gabber whose proof can be found in
  \cite[\S1.1]{Cisinski-Deglise}. The general case follows by the
  Hochschild-Serre spectral sequence.

  For the second we can assume that $X$ is connected; in this case
  \cite[VI.9.8]{AGT} shows that $X_\fet$ is equivalent to the
  classifying topos $B_{\pi_1(X)}$ and that
  \begin{displaymath}
    H^n(X_\fet,F)\simeq H^n(\pi_1(X),F)
  \end{displaymath}
  for all abelian sheaves $F$ on $X_\fet$. It thus suffices to show
  that $H^n(\pi_1(X),F)$ vanishes for $|n|\gg0$  and all continuous
  $\pi_1(X)$-modules $F$.

  Serre \cite[Ch. 1 \S2.2 Cor. 3]{Serre} showed that $H^n(\pi_1(X),F)$
  is torsion for $n>0$, so $H^n(\pi_1(X),F)=0$ if $F$ is a
  $\bQ$-module and $n>0$. Again by Abbes-Gros-Tsuji
  \cite[Prop. VI.9.12]{AGT} $X_\fet$ is a coherent topos and thus
  $H^n(\pi_1(X),\_)$ commutes with filtered inductive limits. A
  devissage using the exact sequence
  \begin{displaymath}
    0\to F_\tor\to F\to F\tens\bQ\to F\tens(\bQ/\bZ)\to 0
  \end{displaymath}
  reduces to the case when $F$ is torsion. Since $\pi_1(X)$ acts
  continuously on $F$, the latter is a filtered inductive limit of
  $\pi_1(X)$-modules corresponding to locally constant constructible
  torsion sheaves. So we can assume that $F$ is constructible, and as
  $X$ is affine we can use a result of Achinger
  \cite[Thm. 1.1.1]{Achinger} showing that there are isomorphisms
  \begin{displaymath}
    H^\cdot(X_\et,F)\simeq H^\cdot(\pi_1(X),F).
  \end{displaymath}
  Thus the assertion follows from the \'etale case.
\end{demo}
In what follows $\Lambda_X$ denotes the sheafification in $\Fet(X)$ of
the constant presheaf with value $\Lambda$ (here $X$ is any
scheme). If $\pi:Y\to X$ is an object of $\Fet(X)$ then
$\Lambda_Y=\pi^*\Lambda_X$. 

\begin{prop}\label{prop:globally-perfect}
  If $M$ is a globally perfect object of $D^+(X_\fet,\Lambda)$ then
  $M$ has finite Tor-dimension. Suppose conversely that $k_\et$ has
  finite cohomological dimension, $X$ is affine and $M$ satisfies
  \begin{itemize}
  \item $M$ has finite Tor-dimension, and
  \item for every $\pi:Y\to X$ in $\Fet(X)$, $R\Gamma(Y,\pi^*M)$ is in
    $D_\perf(\Lambda)$.
  \end{itemize}
  Then $M$ is globally perfect.
\end{prop}
\begin{demo}
  Suppose first that $M$ is globally perfect, and let $a\le b$ be as
  in definition \ref{defn:globally-perfect}. If $\bar x$ runs through
  the set of geometric points of $X$, the family of fiber functors
  associated to the points $\rho_X(\bar x)$ of $X_\fet$ is
  conservative by \cite[Lemma VI.9.6]{AGT}. Since
  $(M\tens_\Lambda N)_x\simeq M_x\tens_\Lambda N$ for any such point
  it suffices to show that $M_x$ has finite Tor-dimension as a
  $\Lambda$-module for all $x$, but this is true since $M_x$ is a
  direct limit of complexes in $D^{[a,b]}_\fTd(\Lambda)$.

  Suppose conversely that $X$ is affine, $k_\et$ has finite
  cohomological dimension and $M$ satisfies the two conditions. We can
  then apply lemma \ref{lemma:coh-dim} to conclude that $X_\fet$ has
  finite cohomological dimension. As in the proof of lemma
  \ref{lemma:coh-dim} we can identify $X_\fet$ with $B_{\pi_1(X)}$ and
  $M$ with an object of $D(\Lambda[\pi_1(X)])$ whose image in
  $D(\Lambda)$ belongs to $D^{[a,b]}_\fTd(\Lambda)$. If $Y\to X$ is in
  $\Fet(X)$, $\pi_1(Y)$ is a subgroup of $\pi_1(X)$ and thus
  \begin{displaymath}
    \cohdim(Y_\fet)\le \cohdim(X_\fet)
  \end{displaymath}
  by \cite[3.3 Prop. 14]{Serre}. Since $X_\fet$ has finite
  cohomological dimension we can apply \cite[Exp. XVII
  Thm. 5.2.11]{SGA4} to conclude that $R\Gamma(Y,M)$ is in
  $D^{[a,b+\cohdim(X)]}_\fTd(\Lambda)$, and since $R\Gamma(Y,\pi^*M)$
  is in $D_\perf(\Lambda)$ by hypothesis, $M$ is globally perfect.
\end{demo}

\begin{lemma}\label{lemma:funct-isom}
  For any $\pi:Y\to X$ in $\Fet(X)$ and $M$ in $D^+(X_\fet,\Lambda)$
  there is a functorial isomorphism
  \begin{displaymath}
    M\tens_{\Lambda_X}^L R\pi_*\Lambda_Y\isom 
    R\pi_*\pi^*(M).
  \end{displaymath}
\end{lemma}
\begin{demo}
  The adjunction $\pi^*R\pi_*(\Lambda_Y)\to\Lambda_Y$ yields a
  morphism
  \begin{align*}
    \pi^*(M\tens_{\Lambda_X}^LR\pi_*(\Lambda_Y))
    &\simeq
      \pi^*(M)\tens_{\Lambda_Y}^L\pi^*R\pi_*(\Lambda_Y)\\
    &\to
      \pi^*(M)\tens_{\Lambda_Y}^L\Lambda_Y\simeq\pi^*(M)
  \end{align*}
  and applying the adjunction to this yields the morphism in the
  lemma. It will be an isomorphism if it is after pulling it back by a
  covering morphism $g:W\to X$ in $\Fet(X)$. By \cite[VI.9.4]{AGT}
  $g^*$ has an exact left adjoint, so $g^*$ sends injectives to
  injectives and it follows that $g^*$ and $R\pi_*$ commute.  We may
  therefore replace $\pi:Y\to X$ by its base change $W\times_XY\to
  W$. Now there is a finite \'etale $g:W\to X$ such that $W\times_XY$
  is a disjoint sum of copies of $W$, and in this case the assertion is
  clear.
\end{demo}

The same argument as in the proof of the last lemma shows:

\begin{lemma}\label{lemma:flatness}
  Suppose $\pi:Y\to X$ is finite \'etale and Galois with group
  $G$. The natural morphism
  \begin{displaymath}
    \pi_*\Lambda_Y\to R\pi_*\Lambda_Y
  \end{displaymath}
  is an isomorphism, and $\pi_*\Lambda_Y$ with its natural $G$-action
  is a flat $\Lambda_X[G]$-module.
\end{lemma}

\noindent \textit{Proof of theorem \ref{thm:perfection}:} By lemma
\ref{lemma:flatness}.  the sheaf
\begin{displaymath}
  \cG:=\pi_*\Lambda_Y
\end{displaymath}
is a flat sheaf of $\Lambda_X[G]$-modules. Since $M$ is globally
perfect it has finite Tor-dimension as an object of
$D(X_\fet,\Lambda_X)$ by proposition \ref{prop:globally-perfect}.
Therefore $M\tens^L_{\Lambda_X}\cG$ has finite Tor-dimension as an
object of $D(X_\fet,\Lambda_X[G])$. Since $X_\fet$ has finite
cohomological dimension by lemma \ref{lemma:coh-dim} we can appeal
once again to \cite[Exp.\ XVII Thm.\ 5.2.11]{SGA4} to conclude that
\begin{displaymath}
  M_Y:=R\Gamma(X_\fet,M\tens^L_{\Lambda_X}\cG)
\end{displaymath}
is in $D_\fTd(\Lambda[G])$. By the same reasoning the image of $M_Y$
in $D(\Lambda)$ lies in $D^b(\Lambda)$; then $M_Y$ is in
$D^b(\Lambda[G])$. Since $\Lambda$ is noetherian it follows that $M_Y$
is an object of $D_\perf(\Lambda[G])$.

On the other hand lemmas \ref{lemma:flatness} and
\ref{lemma:funct-isom} show that
\begin{displaymath}
  M\tens^L_{\Lambda_X}\cG\simeq M\tens^L_{\Lambda_X}R\pi_*{\Lambda_Y}
  \simeq R\pi_*\pi^*M
\end{displaymath}
in $D(X_\fet,\Lambda_X)$, whence isomorphisms
\begin{displaymath}
  M_Y=R\Gamma(X_\fet,M\tens^L_{\Lambda_X}\cG)
  \simeq R\Gamma(X_\fet,R\pi_*\pi^*M)
  \simeq R\Gamma(Y_\fet,\pi^*M).
\end{displaymath}
Since these isomorphisms are functorial and the $G$-module structure
of $\cG=\pi_*\Lambda_Y$ was induced by the action of $G$ on $Y$ the
induced isomorphisms
\begin{displaymath}
  H^\cdot(M_Y)\simeq H^\cdot(Y_\fet,\pi^*M)
\end{displaymath}
are isomorphisms of $\Lambda[G]$-modules.

\section{The main theorem}
\label{sec:main-theorem}

Denote by $K_0$ the fraction field of a Cohen ring of $k$, and recall
that $K$ is the fraction field of $\Lambda$, which we can make into an
extension field of $K_0$. In what follows rigid cohomology has
coefficients in $K_0$.

\begin{lemma}
  Suppose $f:X\to Y$ is a morphism of schemes. If $M$ is a globally
  perfect object of $D(X_\fet,\Lambda)$, $Rf_*M$ is a globally perfect
  object of $D(Y_\fet,\Lambda)$.
\end{lemma}
\begin{demo}
  If $V\to Y$ is in $\Fet(Y)$, $X\times_YV\to X$ is in $\Fet(X)$ and
  \begin{displaymath}
    R\Gamma(V,Rf_*M)\simeq R\Gamma(X\times_YV,M)
  \end{displaymath}
  and the assertion follows.
\end{demo}

\begin{thm}\label{thm:main}
  Suppose $k$ is a field of characteristic $p>0$ and $\Lambda$ is a
  complete noetherian local ring with residue field $k$ and fraction
  field $K$. Suppose $X$ is a smooth affine curve over $k$. Neither
  of the functors
  \begin{displaymath}
    H^\cdot_\rig(\_)\tens_{K_0}K,\ H^\cdot_{\rig,c}(\_)\tens_{K_0}K:
    \Fet(X)\to\Mod_K^\cdot
  \end{displaymath}
  have a globally perfect model.
\end{thm}

\noindent \textit{Proof}: By duality it suffices to show this for
cohomology with compact supports, so we treat the case of
$H^\cdot_{\rig,c}(\_)$. We first reduce to the case when $k_\et$ has
finite cohomological dimension. By standard arguments there is a
subfield $k_0\sset k$ finitely generated over $\bF_p$ and a smooth
affine curve $X_0$ over $k_0$ such that $X\simeq
X_0\tens_{k_0}k$. Suppose $M$ is a globally perfect model of
$H^\cdot_{\rig,c}(\_):\Fet(X)\to\Mod_K^\cdot$; if $f:X\to X_0$ is the
projection then
\begin{displaymath}
  R\Gamma(X_0,Rf_*M)\simeq R\Gamma(X,M)
\end{displaymath}
and $Rf_*M$ is a globally perfect model of
\begin{displaymath}
  H^\cdot_{\rig,c}(\_)\tens_WK:\Fet(X_0)\to\Mod_K^\cdot
\end{displaymath}
where $W$ is a Cohen ring for $k_0$. Thus it suffices to show that
this functor does not have a globally perfect model; i.e.\ we can
replace $X/k$ by $X_0/k_0$. Equivalently we can assume that $k_\et$
has finite cohomological dimension, and as $X$ is affine we can apply
theorem \ref{thm:perfection}.

Any smooth affine curve $X$ has an Artin-Schreier cover $\pi:Y\to X$
that ramifies at infinity. So it suffices to invoke the following
lemma:

\begin{lemma}
  With the hypotheses of the theorem, suppose
  \begin{itemize}
  \item $\pi:Y\to X$ is a finite \'etale Galois cover whose group
    $G$ is a $p$-group;
  \item $M$ is a globally perfect model of the functor
    $H^\cdot_\rig(\_)\tens_{K_0}K:\Fet(X)\to\Mod_K^\cdot$
    or of $H^\cdot_{\rig,c}(\_)\tens_{K_0}K:\Fet(X)\to\Mod_K^\cdot$.
  \end{itemize}
  Then $\pi:Y\to X$ extends to a finite \'etale morphism
  $\bar\pi:\bar Y\to \bar X$ of smooth projective curves.
\end{lemma}
\begin{demo}
  We first consider the case of the functor
  $H^\cdot_{\rig,c}(\_)\tens_{K_0}K$. If it has a globally perfect
  module then by the previous theorem there is a perfect complex
  $M_Y^\cdot$ of $\Lambda[G]$-modules such that
  \begin{displaymath}
    H^n(M_Y^\cdot\tens_\Lambda K)\simeq H^n_{\rig,c}(Y)\tens_{K_0}K
  \end{displaymath}
  and then
  \begin{displaymath}
    H^n(M_Y^\cdot)^G\simeq H^n_{\rig,c}(X)\tens_{K_0}K
  \end{displaymath}
  since $\pi:Y\to X$ is finite \'etale.

  Since $G$ is a $p$-group, the group ring $\Lambda[G]$ is local, so a
  finitely generated projective $\Lambda[G]$-module is free. Therefore
  these last isomorphisms imply that
  \begin{equation}\label{eq:chi-formula}
    \chi_c(X)=|G|\chi_c(Y)
  \end{equation}
  where 
  \begin{displaymath}
    \chi_c(Z)=\sum_i(-1)^i\dim_{K_0}H^i_{\rig,c}(Z)
  \end{displaymath}
  for any separated $k$-scheme $Z$ of finite type. 

  Comparing \ref{eq:chi-formula} with the Grothendieck-Ogg-Shafarevich
  formula shows that the ramification of $\pi:Y\to X$ is tame at
  infinity. But since $G$ is a $p$-group, $\pi$ must be totally wild
  at infinity if it is ramified at all. It is therefore unramified at
  infinity, which is the conclusion of the lemma.

  Finally if $H^\cdot_\rig(\_)\tens_{K_0}K$ has a global perfect model
  the same argument shows that $\chi(Y)=|G|\chi(X)$ where $\chi$ is
  now the Euler characteristic for rigid cohomology without
  supports. But for any smooth variety $Z$, $\chi(Z)=\chi_c(Z)$ and we
  conclude as before.
\end{demo}

\section{Examples}
\label{sec:examples}

The ``standard'' cohomology theories can all be computed by complexes
of \'etale sheaves, and therefore satisfy \'etale cohomological
descent. Thus when finite they have globally perfect models. For the
sake of completeness let's review the basic examples. The reader
should have no problem thinking of global perfect models for various
sorts of cohomology theories: $\ell$-adic cohomology with or without
supports, crystalline cohomology of proper smooth varieties (use the
crystalline-\'etale topos $(X/V)_\cryset$ of Berthelot-Breen-Messing
\cite{BBM}, or else the de Rham-Witt complex). The case of $p$-adic
\'etale cohomology with compact supports (where $p$ is the residual
characteristic) was used in \cite{Crew}

Suppose on the other hand that $k$ is perfect and
$W^\dagger\Omega^\cdot_{X/W}$ is the complex of Davis, Langer and Zink
\cite{DLZ} and let $\Lambda=W(k):=W$. If the $W$-modules
$H^\cdot(Y,W^\dagger\Omega^\cdot_{Y/W})$ were finitely generated for all
$Y\to X$ in $\Fet(X)$ then $W^\dagger\Omega^\cdot_{Y/W}$ would be a
global perfect model of $H^\cdot_\rig(\_):\Fet(X)\to\Mod^\cdot_K$
where $K$ is the fraction field of $W=W(k)$. We conclude that for
every smooth affine curve $X$ there is a finite \'etale $Y\to X$ such
that $H^\cdot(Y,W^\dagger\Omega^\cdot_{Y/W})$ is not a finitely
generated $W$-module.

\bigskip

\parindent=0pt

Tomoyuki Abe:\\
Kavli Institute for the Physics and Mathematics of the Universe (WPI)\\
The University of Tokyo,
5-1-5 Kashiwanoha, Kashiwa, Chiba, 277-8583, Japan\\
e-mail: \texttt{tomoyuki.abe@ipmu.jp}
\bigskip

Richard Crew:\\
Department of Mathematics, University of Florida\\
358 Little Hall, Gainesville FL 32611\\
USA\\
e-mail: \texttt{rcrew@ufl.edu}

\end{document}